\documentclass[12pt]{article}
\usepackage[utf8]{inputenc}
\usepackage{amsmath}
\usepackage{amssymb}
\usepackage{mathrsfs}
\usepackage{amsthm}
\usepackage{cite}
\usepackage{hyperref}
\usepackage{enumerate}
\usepackage[left=0.7in,right=0.7in,top=0.7in,bottom=0.7in]{geometry}
\usepackage{titlesec}
\titleformat*{\section}{\large\bfseries}
\newtheorem{theorem}{Theorem}[section]
\newtheorem{lemma}[theorem]{Lemma}
\newtheorem{corollary}[theorem]{Corollary}
\newtheorem{definition}[theorem]{Definition}
\newtheorem{example}[theorem]{Example}

\numberwithin{equation}{section}
\title{Isometries of the product of composition\\ operators on weighted Bergman space}
\author{\large Anuradha Gupta and Geeta Yadav$^*$}
\date{}
\begin{document}
\maketitle
\begin{abstract} 
In this paper the necessary and sufficient conditions for the product of composition operators to be isometry are obtained on weighted Bergman space. With the help of a counter example we also proved that unlike on $\mathcal{H}^2(\mathbb{D})$ and $\mathcal{A}_{\alpha}^2(\mathbb{D}),$ the composition operator on $\mathcal{S}^2(\mathbb{D})$ induced by an analytic self map on $\mathbb{D}$ with fixed origin need not be of norm one. We have generalized the Schwartz's \cite{Schwartz Thesis} well known result on $\mathcal{A}_{\alpha}^2(\mathbb{D})$ which characterizes the almost multiplicative operator on $\mathcal{H}^2(\mathbb{D}).$ 


\end{abstract}

\section{Introduction and Preliminaries}

Let $\mathbb{D}=\{z\in \mathbb{C} :\, |z|<1\}$ denote the open unit disk in the complex plane $\mathbb{C}$ and $H(\mathbb{D})$ denote the space of all analytic complex valued functions on $\mathbb{D}$. Let $\mathcal{H}^2(\mathbb{D}):=\mathcal{H}^2$ denote the Hardy space defined by
$$\mathcal{H}^2(\mathbb{D})=\Big\{ f(z)=\sum_{n=0}^{\infty} a_n z^n  \in H(\mathbb{D}): \,||f||_{\mathcal{H}^2}^2= \sum_{n=0}^{\infty} |a_n|^2 < \infty\Big\},$$
the derivative Hardy space  denoted by $\mathcal{S}^2(\mathbb{D}):=\mathcal{S}^2$ is defined by
\begin{align*}
 \mathcal{S}^2(\mathbb{D})&=\Big\{f\in H(\mathbb{D}): \, f'\in \mathcal{H}^2(\mathbb{D}) \Big\}\\
                &=\Big\{f(z)=\sum_{n=0}^{\infty} a_n z^n  \in H(\mathbb{D}): \,||f||_{S^2}^2= |a_{0}|^2+\sum_{n=1}^{\infty} n^2|a_n|^2 < \infty \Big\}.
\end{align*}
	
For $\alpha>-1,$ let $\mathcal{A}_{\alpha}^2(\mathbb{D}):=\mathcal{A}_{\alpha}^2$ denote the weighted Bergman space defined by 
$$ \mathcal{A}_{\alpha}^2(\mathbb{D})=\{f \in H(\mathbb{D}):\, ||f||_{\alpha}^2=\int_{\mathbb{D}} |f(z)|^2 \, dA_{\alpha}(z) \,< \infty \},$$ 
where $dA_{\alpha}(z)=(\alpha+1)(1-|z|^2)^{\alpha}\, dA(z)$ and $dA$ denote the normalized  area measure on $\mathbb{D}.$ 
For $\alpha=0,$ $\mathcal{A}_{\alpha}^2(\mathbb{D})$ is called the Bergman space and it is denoted by $\mathcal{A}^2(\mathbb{D}):=\mathcal{A}^2.$
It is well known that $e_n(z)=\sqrt{\frac{\Gamma(n+2+\alpha)}{n! \, \Gamma(2+\alpha)}}\, z^n$ for $n=0,1,...,$ where $\Gamma$ is gamma function, forms an orthonormal basis for $\mathcal{A}_{\alpha}^2.$ Clearly, $e_0=1.$  By $\mathbb{Z_+},$ we denote the set of non negative integers. \\

The reproducing kernel denoted by $K_w \in \mathcal{A}_{\alpha}^2$ for point evaluation at $w$ in $\mathbb{D}$ which satisfies $f(w)=\langle f, K_w \rangle$  for each $f \in \mathcal{A}_{\alpha}^2(\mathbb{D})$ is 

\begin{equation} \label{Kernel of Bergman}
 K_w(z)=\frac{1}{(1-\overline{w}z)^{(\alpha+2)}}
\end{equation}
and 
\begin{equation} \label{Norm of Kernel of Bergman}
||K_w||^2_{\alpha}=\frac{1}{(1-|w|^2)^{(\alpha+2)}}.
\end{equation}

By $Aut(\mathbb{D})$ we denote the set of all disc automorphism, that is, the set of all maps in $H(\mathbb{D})$ which maps $\mathbb{D}$ into $\mathbb{D}$ and are bijective. Equivalently, $\Phi\in Aut(\mathbb{D})$ if and only if there exist $b\in \mathbb{D},$ $\gamma \in \mathbb C$  such that $|\gamma|=1$ and
$$ \Phi(w)= \gamma \frac{b-w}{1-\bar{b}w},\,\, w\in \mathbb{D}  $$

Let $\mathcal{H}^\infty(\mathbb{D}):=\mathcal{H}^\infty$ denote the space of all bounded analytic functions on the open unit disk with the supremum norm.\\

Let $\Phi\in H(\mathbb{D})$ and  $\Phi$ be a self map on $\mathbb{D},$ that is, $\Phi:\mathbb{D} \longrightarrow \mathbb{D}.$ Then the composition operator induced by $\Phi$ on $H(\mathbb{D})$ is defined as
$$ C_{\Phi}f=f\circ \Phi\,\, \text{for}\, f\in H(\mathbb{D}) .$$\
Let $\Phi, u \in H(\mathbb{D})$ and  $\Phi$ be a self map on $\mathbb{D}.$  Then the weighted composition operator on $H(\mathbb{D})$ is defined as
$$ C_{u,\Phi}f=u\cdot (f\circ \Phi)\,\, \text{for}\, f\in H(\mathbb{D}) .$$
If $u \equiv 1$ on $\mathbb{D}$ then the weighted composition operator is the composition operator $C_{\Phi}$ induced by $\Phi.$
\\

Throughout the paper we assume that $\Phi, \Psi \in H(\mathbb{D})$ and are self maps on $\mathbb{D}$.
Carswell and Hammond gave the characterization for isometric composition operator on weighted Bergman space. From (Propostion 3.1, \cite{Maximal norm}) it follows that $C_\Phi:\mathcal{A}_ {\alpha}^{2} \longrightarrow \mathcal{A}_ {\alpha}^{2}$  is an isometry if and only if $\Phi$ is a rotation. One can easily prove that if $\Phi:\mathbb{D} \longrightarrow \mathbb{D}$ is a rotation, that is, there exists $\lambda\in \mathbb{C}$ with $|\lambda|=1$ such that  $\Phi(z)=\lambda z$ for $z\in \mathbb{D},$ then  $C_{\Phi}^*=C_\Psi$ where $\Psi(z)=\overline{\lambda}z$ for $z\in \mathbb{D}.$ In fact, $C_{\Phi}^*C_{\Phi}=I$ and $C_{\Phi}C_{\Phi}^*=I$ where $I$ is identity operator on $\mathcal{A}_{\alpha}^2.$ Consequently, $C_{\Phi}^{-1}=C_{\Phi}^*$, $C_{\Phi}$ is unitary and we have the following result:

\begin{theorem} \label{CH Bergman}
A composition operator $C_\Phi$ on $\mathcal{A}_{\alpha}^2$ is an isometry if and only if it is unitary. 
\end{theorem} 
A function $g$ with domain $X$ is said to be bounded away from zero on $X$ if there exists a positive constant $c$ which satisfy $|g(z)|\geq c$ for each $z$ in domain $X.$ \\

In \cite{Fredholmness} Shaabani gave the following characterization for invertibility of weighted composition operator on the weighted Bergman space which is instrumental in proving the results in the subsequent sections. 
\begin{theorem}  

Let $\Psi$ be an analytic map on $\mathbb{D}$, and let $\Phi$ be an analytic self-map on $\mathbb{D}$. The weighted composition operator $C_{\Psi,\Phi}$ is invertible on $\mathcal{H}^2$ or $\mathcal{A}_{\alpha}^2$ if and only if $\Phi \in Aut(\mathbb{D})$ and $\Psi \in \mathcal{H}^\infty$  is bounded away from zero on $\mathbb{D}$. 
\end{theorem}

One can easily note that if $\Psi \equiv 1$  on $\mathbb{D}$ then the composition operator $C_{\Phi}$ is invertible on $\mathcal{A}_{\alpha}^2$ if and only if $\Phi \in Aut(\mathbb{D}).$

Recently, Shaabani \cite{Shaabani Closed} discussed the necessary and sufficient condition for the product $C_{u_{1},\Phi_{1}}C_{u_{2},\Phi_{2}}^*$ and $C_{u_{1},\Phi_{1}}^*C_{u_{2},\Phi_{2}}$ to be invertible on $\mathcal{A}_{\alpha}^2.$ In section 2, we characterize the isometry operators $C_\Phi C_\Psi^*$ and discuss the necessary condition for operators $C_\Phi^* C_\Psi$ to be isometry. With the help of a counter example we proved that unlike $\mathcal{H}^2$ and $\mathcal{A}_{\alpha}^2$ if an analytic self map $\Phi$ on $\mathbb{D}$ with $\Phi(0)=0$ induces the composition operator $C_{\Phi}$ on  $\mathcal{S}^2$ then the norm of $C_{\Phi}$ may not be one. In section 3, we generalizes the Schwartz's  \cite{Schwartz Thesis}  well known result on $\mathcal{A}_{\alpha}^2$  which characterizes the almost multiplicative operators on $\mathcal{H}^2.$ This result helps us in characterizing the isometry operators $C_\Phi^* C_\Psi$.

\section{The products $C_{\Phi}C_{\Psi}^*$ and $C_{\Phi}^*C_{\Psi}$}

We begin this section with the definition of adjoint of composition operator on $\mathcal{A}_{\alpha}^2$  which will be used repeatedly. 
\begin{definition}
The adjoint operator of a composition operator $C_\Phi$ on $\mathcal{A}_{\alpha}^2$ denoted by $C_\Phi^*$ is the operator 
which satisfy
$$\langle C_{\Phi}f, g \rangle=\langle f,C_{\Phi}^*g \rangle $$
for $f,g \in \mathcal{A}_{\alpha}^2.$ Also,
\end{definition}
$$\langle f,C_{\Phi}^* K_w  \rangle= \langle C_{\Phi}f, K_w  \rangle= (f\circ\Phi)(w)=f(\Phi(w))=\langle f,K_{\Phi(w)}  \rangle$$
for $f\in \mathcal{A}_{\alpha}^2.$ Thus, 
\begin{equation} \label{adjoint property of Comp Op}
C_{\Phi}^* K_w=K_{\Phi(w)}
\end{equation}

Next we will generalize the Carswell and Hammond result for $C_{\Phi}C_{\Psi}^*$ to be an isometry.\\

\begin{theorem} \label{Com Op adjoint Comp Op isometry iff}
The product $C_{\Phi}C_{\Psi}^*$ is an isometry if and only if both the composition operators $C_\Phi$ and $C_\Psi$ are isometries, equivalently, the inducing maps $\Phi$ and $\Psi$ are rotations.
\end{theorem}
\textbf{Proof.} Let $C_{\Phi}C_{\Psi}^*$ be an isometry. From the equivalent conditions of isometry on Hilbert space it follows that  
\begin{equation} \label{TS* isometry theorem}
C_{\Psi}(C_{\Phi}^*C_{\Phi}C_{\Psi}^*)=(C_{\Phi}C_{\Psi}^*)^*(C_{\Phi}C_{\Psi}^*)=I. 
\end{equation}
 Thus, $C_{\Psi}$ is surjective. Hence $C_{\Psi}$ is injective because $\Psi$ is non constant. Therefore $C_{\Psi}$ is invertible and $\Psi$ is disc automorphism. So there exists $w\in \mathbb{D}$ such that $\Psi(w)=0.$ From Equation (\ref{Kernel of Bergman}) $K_{0}=1$ on $\mathbb{D}.$ Clearly, $C_{\Phi}K_{0}=K_{0}.$ From Equations (\ref{Norm of Kernel of Bergman}) and (\ref{adjoint property of Comp Op}) we get
\begin{equation*}
||C_{\Phi}C_{\Psi}^* K_w||_{\alpha}^2=||C_{\Phi}K_0||_{\alpha}^2=1. 
\end{equation*}
Since $C_{\Phi}C_{\Psi}^*$ is an isometry we have
\begin{equation*}
1=||C_{\Phi}C_{\Psi}^* K_w||_{\alpha}^2=||K_w||_{\alpha}^2=\frac{1}{(1-|w|^2)^{(\alpha+2)}}.
\end{equation*}
Thus, $w=0.$ This implies $\Psi$ is a rotation because it is a disc automorphism and $\Psi(0)=0.$ Hence $C_{\Psi}$ is an isometry and unitary as well. Now by pre and post multiplying Equation (\ref{TS* isometry theorem}) with $C_{\Psi}^*$ and $C_{\Psi}$ respectively we get $C_{\Phi}^*C_{\Phi}=I,$ that is $C_{\Phi}$ is an isometry. 
          
For the converse part, let $C_\Phi$ and $C_{\Psi}$ be isometries. Then both $C_\Phi$ and $C_{\Psi}^*$ are unitary. Since the product of two unitary operators is unitary we get $C_{\Phi}C_{\Psi}^*$ is unitary and hence isometry.  \\

The following result is direct consequence of Theorem \ref{Com Op adjoint Comp Op isometry iff}.
\begin{corollary}
The product $C_{\Phi}C_{\Psi}^*$ is unitary if and only if $C_\Phi$ and $C_\Psi$ are isometry.
\end{corollary}

\begin{theorem} \label{theorem Product Com Op norm one iff induced map fix origin}
The product $C_{\Phi}^*C_{\Psi}$ has norm one if and only if $\Phi(0)=0=\Psi(0).$ 
\end{theorem} 
\textbf{Proof.} Let $||C_{\Phi}^*C_{\Psi}||=1.$ Clearly, applying Equation (\ref{adjoint property of Comp Op}) with the fact that $C_{\Psi}K_{0}=K_{0}$ we get  
\begin{equation*}
C_{\Phi}^*C_{\Psi}K_{0}=C_{\Phi}^*K_{0}=K_{\Phi(0)}.
\end{equation*}
Then 
$$\frac{1}{(1-|\Phi(0)|^2)^{(\alpha+2)}}=||C_{\Phi}^*C_{\Psi}K_{0}||_{\alpha}\leq ||C_{\Phi}^*C_{\Psi}|| \, ||K_{0}||_{\alpha}=1$$
so $\Phi(0)=0$ and $||C_{\Psi}^*C_{\Phi}||=||C_{\Phi}^*C_{\Psi}||=1.$  Similarly, $||C_{\Psi}^*C_{\Phi}||=1$ gives $\Psi(0)=0.$\\

For the converse part, let $\Phi(0)=0$ and $\Psi(0)=0.$ Then using the upper bound on the norm of a composition operator by (Theorem $11.6$, \cite{K Zhu Op Theory in fun spaces})(with $p=2$) we have 
\begin{equation*}
||C_{\Phi}^*C_{\Psi}|| \leq ||C_{\Phi}||\, ||C_{\Psi}|| \leq \left(\frac{1+|\Phi(0)|}{1-|\Phi(0)|}\right)^{(\alpha+2)/2}\left(\frac{1+|\Psi(0)|}{1-|\Psi(0)|}\right)^{(\alpha+2)/2}=1.
\end{equation*}
For reverse inequality, consider
$$ ||C_{\Phi}^*C_{\Psi}||\geq \frac{||C_{\Phi}^*C_{\Psi}K_0||_{\alpha}}{||K_0||_{\alpha}}=\frac{||K_{\Phi(0)}||_{\alpha}}{||K_0||_{\alpha}}=\frac{||K_0||_{\alpha}}{||K_0||_{\alpha}}=1.$$  \\

As it is well known that the norm of an isometry operator on any Hilbert space is one but the converse is not true so the above theorem leads to the following result:
\begin{corollary} \label{Compadjoint Comp isom then at zero is zero}
If the product $C_{\Phi}^*C_{\Psi}$ is an isometry then $\Phi(0)=0$ and $\Psi(0)=0.$
\end{corollary}

Unlike $\mathcal{H}^2$ and $\mathcal{A}_{\alpha}^2$, if an analytic self map $\Phi$ on $\mathbb{D}$ that induces the composition operator on $\mathcal{S}^2$ and fixes the origin, then this does not imply that the norm of the composition operator induced by $\Phi$ on $\mathcal{S}^2$  is one. This can be verified from the following example.

\begin{example} \label{Example Cphi Com op}
Let $\Psi$ be an analytic, self map defined by $\Psi(z)=z^3\,\, \forall\, z\in \mathbb{D}$ so that $\Psi(0)=0.$ Clearly, $\Psi' \in \mathcal{H}^\infty$ and this implies that $C_\Psi$ is a composition operator on $\mathcal{S}^2$ (see \cite{Roan Paper}).  
Define $f(z)=z,\,\, z\in \mathbb{D}$ so that $f\in \mathcal{S}^2.$\\
Then 
\begin{align*}
||C_\Psi|| &\ge\frac{||C_\Psi f||_{\mathcal{S}^2}} {||f||_{\mathcal{S}^2}}\\     
                            &= \frac{||\Psi||_{S^2}}{||f||_{\mathcal{S}^2}}=\frac{3}{1}=3
\end{align*}
Thus, we have proved that $\Psi(0)=0$ but $||C_\Psi||\ne 1.$ 
\end{example}

The above example is instrumental in showing that unlike $\mathcal{A}_{\alpha}^2$ the converse part of the statement in Theorem \ref{theorem Product Com Op norm one iff induced map fix origin} does not hold on $\mathcal{S}^2$ in the following example. 
	
\begin{example}
Let $\Psi$ be as defined in the Example \ref{Example Cphi Com op} and $\Phi$ be defined by $\Phi(z)=z \,\, \forall\,  z\in \mathbb{D}.$ Then $\Phi$ is analytic, self map on $\mathbb{D}$ with $\Phi(0)=0$ and $\Phi' \in \mathcal{H}^\infty.$ So $C_\Phi$ and $C_\Psi$ are composition operators on $\mathcal{S}^2$ (see \cite{Roan Paper}). Clearly, $C_\Phi f=f$ for $f\in \mathcal{S}^2$ and this gives 
$$\langle C_\Phi f, g\rangle=\langle f  , g\rangle=\langle f  , C_\Phi g\rangle \,\, \text{for}\, f,g\in \,\mathcal{S}^2 .$$
 Thus, $C_\Phi$ is an isometry and $C_\Phi^*=C_\Phi.$ Consequently, for $f\in \mathcal{S}^2$ 
\begin{equation*} \label{equation in example of converse part them}
||C_\Phi^*C_\Psi f||_{\mathcal{S}^2}=||C_\Psi f||_{\mathcal{S}^2}
\end{equation*}
Using Example \ref{Example Cphi Com op} we get
\begin{align*}
||C_\Phi^*C_\Psi||&=Sup_{f\ne 0\in \mathcal{S}^2}\frac{||C_\Phi^*C_\Psi f||_{\mathcal{S}^2}} {||f||_{\mathcal{S}^2}}\\
          &=Sup_{f\ne 0\in \mathcal{S}^2}\frac{||C_\Psi f||_{\mathcal{S}^2}} {||f||_{\mathcal{S}^2}}\\
					&=||C_\Psi||\\
					&\geq 3
\end{align*}
Thus we have proved that $\Phi(0)=0$ and $\Psi(0)=0$ but $||C_\Phi^*C_\Psi||\neq 1.$ 
\end{example}

\section{Almost Multiplicative}

In this section we will discuss the necessary and sufficient condition for operators $C_{\Phi}^*C_{\Psi}$ to be isometry on  $\mathcal{A}_{\alpha}^2.$ 
\begin{definition}
An operator $T$  on a Hilbert space $\mathcal{H}$ is said to be almost multiplicative if 
$$T(f_1 \cdot f_2)= Tf_1 \cdot Tf_2$$ 
whenever $f_1, f_2 \in \mathcal{H}$ satisfy $f_1 \cdot f_2 \in \mathcal{H}$
\end{definition}
 
For our convenience we will write $f_1 f_2$ to denote $ f_1 \cdot f_2.$ In the following result on the weighted Bergman space we will generalize the Schwartz's \cite{Schwartz Thesis} result which shows the relation between almost multiplicative and composition operator on the Hardy space.

\begin{theorem} \label{composition iff almost multiplicative}
Let  $T$ be a non zero bounded linear operator on $\mathcal{A}_{\alpha}^2.$ Then $T$ is almost multiplicative if and only if $T$ is a composition operator i.e, $T=C_{\Phi}$ for some analytic function $\Phi:\mathbb{D} \longrightarrow \mathbb{D}.$ 
\end{theorem}
\textbf{Proof.} Let $T$ be almost multiplicative. First, we will prove that 
$$Te_n=\sqrt{\frac{\Gamma(n+2+\alpha)}{n! \, \Gamma(2+\alpha)}} \, (Tz)^n\,\, \text{for}\,\, n\in \mathbb{Z_+}.$$ 
As $T$ is non zero operator and for $n\in \mathbb{Z_+}$ 
$$Te_n=T(e_n e_0)=(Te_n)(Te_0)$$
so $Te_0$ can not be a zero function. Since
$$Te_0= T(e_0 e_0)=(Te_0) (Te_0)$$
therefore, $Te_0=e_0.$ Now for $n\geq 1$
\begin{align} 
Te_n&=T\left(\sqrt{\frac{\Gamma(n+2+\alpha)}{n! \, \Gamma(2+\alpha)}}\, z^n \right) \notag\\
      &=\sqrt{\frac{\Gamma(n+2+\alpha)}{n! \, \Gamma(2+\alpha)}}\, (Tz)^n \notag\\
			 &=\sqrt{\frac{\Gamma(n+2+\alpha)}{n! \, \Gamma(2+\alpha)}}\, \Phi^n \,\,\text{where}\,\, \Phi=Tz \label{Te_n in terms of phi}
\end{align}
This gives $\Phi^n=\sqrt{\frac{n! \, \Gamma(2+\alpha)}{\Gamma(n+2+\alpha)}}(Te_n)$ for $n\in \mathbb{Z_+}$. Further
\begin{align*}
||\Phi^n||_{\alpha}&=\sqrt{\frac{n! \, \Gamma(2+\alpha)}{\Gamma(n+2+\alpha)}}||(Te_n)||_{\alpha}\\
				   &\leq \sqrt{\frac{n! \, \Gamma(2+\alpha)}{\Gamma(n+2+\alpha)}} ||T||\, ||e_n||_{\alpha}\\
           &=\sqrt{\frac{n! \, \Gamma(2+\alpha)}{\Gamma(n+2+\alpha)}} ||T||.                      
\end{align*}
Since $\sqrt{\frac{n! \, \Gamma(2+\alpha)}{\Gamma(n+2+\alpha)}} \leq 1$ for $n\in \mathbb{Z_+},$ therefore, for $n\in \mathbb{Z_+}$  
\begin{equation*} \label{phi func bounded by T}
||\Phi^n||_{\alpha}\leq ||T||
\end{equation*}
Thus $\{||\Phi^n||_{\alpha}\}$ is a bounded sequence. Since $\Phi=Tz\in \mathcal{A}_{\alpha}^2$ and any analytic self map on $\mathbb{D}$ induces a composition operator on $\mathcal{A}_{\alpha}^2,$ therefore, it is sufficient to prove that $\Phi$ is a self map on $\mathbb{D}.$ Now, we will prove that $|\Phi(z)| < 1$ $\forall\,z\in \mathbb{D}.$\\ 
Let $\delta>0$ be any real number and let $E=\{z\in \mathbb{D}: |\Phi(z)|\geq 1+\delta\}.$ Consider
\begin{align*}
||\Phi^n||_{\alpha}^2 &=\int_{\mathbb{D}} |\Phi(z)|^{2n}\, dA_{\alpha}(z) \\
						 &\geq \int_{E} |\Phi(z)|^{2n}\, dA_{\alpha}(z) \\
						 &\geq \int_{E} (1+\delta)^{2n}\, dA_{\alpha}(z) \\
						&=(1+\delta)^{2n} A_{\alpha}(E).		
\end{align*}   
Since $1+\delta>1,$ from above it follows that measure of set $E$ is zero, i.e, $A_{\alpha}(E)= 0$ otherwise 
$$\lim_{n \to \infty} ||\Phi^n||_{\alpha}=\infty$$
and this is a contradiction to the fact that $\{||\Phi^n||_{\alpha}\}$ is a bounded sequence. So, $|\Phi(z)|\leq 1$ a.e on $\mathbb{D}.$ By combining (Corollary $1.29$, \cite{K Zhu Space of Hol} for $n=1$) with (Corollary $1.30$, \cite{K Zhu Space of Hol} for $p=2$, $n=1$) and applying for analytic function $\Phi$ on $\mathbb{D}$ we get that for  all $b\in \mathbb D$  and $0\leq r<1-|b|$ 
\begin{align*}
|\Phi(b)|^2 &\leq \int _{\mathbb D}\, |\Phi(b+rz)|^2\, dA_{\alpha}(z)\\
          &\leq \int _{\mathbb D}\, \, dA_{\alpha}(z) \\
					&=1
\end{align*}
Thus, $|\Phi(z)|\leq 1 \,\, \forall \,z\in \mathbb{D}.$ Now prove that $|\Phi(z)| \ne 1$ for any $z\in \mathbb{D}.$  

We will prove this by contradiction. 
Let $|\Phi(z')|= 1$ for some $z'\in\mathbb{D}$ then by the Maximum Modulus Principle there exists some $\beta$ with $|\beta|=1$ such that $\Phi(z)=\beta$ for all $z\in \mathbb{D},$ that is, $\Phi$ is a constant  function.  Now using this along with Equation (\ref{Te_n in terms of phi}) we have

$\langle T^* e_{0}, e_n \rangle=\langle  e_{0}, Te_n\rangle=\big \langle  e_{0}, \sqrt{\frac{\Gamma(n+2+\alpha)}{n! \, \Gamma(2+\alpha)}}\, \, \beta^n \big\rangle=\sqrt{\frac{\Gamma(n+2+\alpha)}{n! \, \Gamma(2+\alpha)}}\,\, \bar{\beta}^n$\\
Since $\{e_n\}_{n=0}^{\infty}$ is an orthonormal basis of $\mathcal{A}_{\alpha}^2,$ we can write
\begin{align*} 
T^*e_{0}&=\sum_{n=0}^{\infty} c_n\, e_n\,\text{for some scalars}\,  c_{n}\\
         &=\sum_{n=0}^{\infty} \langle T^* e_{0}, e_n\rangle\, e_n \,\text{where} \, c_{n}=\langle T^* e_{0}, e_n\rangle
\end{align*}

Now Consider
\begin{align*}
||T^*e_{0}||_{\alpha}^2&= \langle T^* e_{0}, T^*e_0\rangle\\
           &=\Big\langle \sum_{n=0}^{\infty} \langle T^* e_{0}, e_n\rangle\, e_n,\sum_{n=0}^{\infty} \langle T^* e_{0}, e_n\rangle\, e_n \Big\rangle\\
							&=\sum_{n=0}^{\infty} \langle T^* e_{0}, e_n\rangle \overline{\langle T^* e_{0}, e_n\rangle}\\
							&=\sum_{n=0}^{\infty} |\langle T^* e_{0}, e_n\rangle|^2\\
							&=\sum_{n=0}^{\infty} \frac{\Gamma(n+2+\alpha)}{n! \, \Gamma(2+\alpha)}|\beta|^{2n}\\
							&=\infty
 \end{align*}
This contradicts the fact that the norm of $T^*e_0 \in \mathcal{A}_{\alpha}^2$ is finite. Hence $|\Phi(z)|< 1,$ $z\in \mathbb{D}.$ From equation (\ref {Te_n in terms of phi}) we have

$\,\,\,\,\,\,\,\,\,\,\,\,\,\,\,\,\,\,\,\,\,\,\,\,\,\,\,\,\,\,~~~~~~~~~~~Te_n=\sqrt{\frac{\Gamma(n+2+\alpha)}{n! \, \Gamma(2+\alpha)}}\, \Phi^n =e_n\circ\Phi= C_{\Phi} e_n$\\
and it gives by linearity and continuity that $T=C_\Phi.$ Hence, we have proved that $T$ is a composition operator. \\
Conversely, let $T$ be a composition operator, i.e, $T=C_{\Phi}.$ Then for $f_1, f_2 \in \mathcal{A}_{\alpha}^2 $ such that $f_1 f_2 \in \mathcal{A}_{\alpha}^2$ and for $z\in \mathbb{D}$ we have
 
$$(C_{\Phi}(f_1f_2))(z)= ((f_1 f_2) \circ \Phi)(z)=(f_1 f_2)(\Phi(z))=f_1(\Phi(z)) f_2(\Phi(z))=((C_{\Phi}f_1) (C_{\Phi}f_2))(z).$$
Thus, $C_{\Phi}(f_1 f_2)=(C_{\Phi}f_1) (C_{\Phi}f_2)$ and $T$ is almost multiplicative.\\


The following result is motivated by Clifford et al. (Lemma 4.4, \cite{Main Paper}) and the proof  follows on the similar lines using Theorem \ref{composition iff almost multiplicative} .
\begin{lemma} \label{Lemma Cpsi=CphiT}
Let $C_{\Psi}=C_{\Phi}S,$ $\Psi$ be a non constant function and let $S$ be a bounded linear operator on $\mathcal{A}_{\alpha}^2.$ Then $S=C_{\beta}$ for some $\beta \in H(\mathbb{D})$ which is a self map on $\mathbb{D}$ and $\Psi=\beta \circ \Phi.$
\end{lemma}

\begin{theorem} \label{Product of adjoint of Comp and adjoint op isometry iff rotation}
The product $C_{\Phi}^*C_{\Psi}$ of operators $C_{\Phi}^*$ and $C_{\Psi}$ is an isometry if and only if $\Phi$ and $\Psi$ are rotations on $\mathbb{D}.$ 
\end{theorem}
\textbf{Proof.} Let $C_{\Phi}^*C_{\Psi}$ be an isometry. It follows from Corollary \ref{Compadjoint Comp isom then at zero is zero} that $\Phi(0)=0$ and $\Psi(0)=0$. Thus, $||C_{\Phi}|| \leq 1$ and $||C_{\Psi}|| \leq 1$. Then, by (Proposition 4.3, \cite{Main Paper})  it follows that $C_{\Psi}$ is an isometry and $C_{\Psi}=C_{\Phi}C_{\Phi}^*C_{\Psi}.$ Since $C_{\Psi}$ is an isometry, therefore $\Psi$ is a rotation and hence, non constant. Applying Lemma \ref{Lemma Cpsi=CphiT} with $S$ as $C_{\Phi}^*C_{\Psi}$ we get that $C_{\Phi}^*C_{\Psi}=C_\beta$ for some $\beta \in H(\mathbb{D})$ which is a self map on $\mathbb{D}$ with $\Psi=\beta \circ \Phi.$ Combining the fact that $C_{\Phi}^*C_{\Psi}$ is an isometry and $C_{\Phi}^*C_{\Psi}=C_\beta$ we get that $\beta$ is a rotation. Since $\Psi$ and  $\beta$ are rotations on $\mathbb{D}$ and $\Psi=\beta \circ \Phi,$ therefore, $\Phi$ is a rotation.\\
Conversely, let $\Phi$ and $\Psi$  are rotations on $\mathbb{D}.$  Then $C_{\Phi}$ and $C_{\Psi}$ are isometry and unitary as well. This implies that $C_{\Phi}^*C_{\Psi}$ is unitary and hence an isometry. \\

The following result is a direct consequence of Theorem \ref{Product of adjoint of Comp and adjoint op isometry iff rotation}.

\begin{corollary}
The product $C_{\Phi}^*C_{\Psi}$ of operators $C_{\Phi}^*$ and $C_{\Psi}$ is unitary if and only if $\Phi$ and $\Psi$ are rotations on $\mathbb{D}.$ 
\end{corollary}

\textbf{Anuradha Gupta}\\
 Department of Mathematics, Delhi College of Arts and Commerce,\\
  University of Delhi, New Delhi-110023, India.\\
  \vspace{0.2cm}
 email: dishna2@yahoo.in\\
   \textbf{Geeta Yadav}\\
  Department of Mathematics, University of Delhi, \\
  New Delhi-110007, India.\\
  email: ageetayadav@gmail.com
\end{document}